\documentclass[12pt]{article}
\usepackage[cp1251]{inputenc}
\usepackage[russian,english]{babel}
\usepackage{latexsym,amsfonts,amssymb,amsmath,amsfonts}
\usepackage{setspace}

 \newtheorem{theorem}{Theorem}[section]

 \numberwithin{equation}{section}

\begin{document}


\begin{center}

{\bf  Inverse problems of recovering lower-order coefficients  from boundary integral data}

\medskip

S.G. Pyatkov, O.A. Soldatov \footnote{ MSC Primary 35R30; Secondary 35K20; 80A20 }

\medskip

\bigskip

\parbox{10cm}{\small
 Abstract.  Under consideration are mathematical models of heat and mass transfer.
We study  inverse problems of recovering lower-order  coefficients in  a second order parabolic equation.
The coefficients  are representable in the form  of  a finite segments of the series
  whose coefficients depending on time are to be determined. The linear case is also considered.
  The overdetermination conditions are the integrals over the boundary of the domain of a solution with weights. 
The main attention is paid to existence,  uniqueness, and stability estimates for
solutions to inverse problems of this type. The problem is reduced to an operator equation which
is studied with the use of the fixed point theorem and a priori estimates.   A solution has all generalized derivatives occurring into the equation summable to some power.  
The method of the proof is constructive and it can be used for developing new numerical algorithms of solving the problem.
 }

\end{center}



\section{Introduction}

We examine the question on identification of lower-order coefficients and the right-hand side in a parabolic equation.  
An equation is of the form 
\begin{equation}\label{e1}
	Mu= u_t+Au=f(t,x), \ \ (t,x)\in Q = (0,T) \times G,
\end{equation}
where   $G\in\mathbb R^n$ is a bounded domain with boundary  $\Gamma$. Let   $S= (0,T)\times\Gamma$.
The function  $f$ and a second order  elliptic operator   $A$ are representable as 
\begin{multline*}
	A(t,x,D)=A_{0}(t,x,D_x)+ \sum_{i=1}^{r}q_i(t)A_i(t,x,D_x) ,\   f=f_0(t,x)+ \sum_{i=r+1}^{s}f_i(t,x)q_i (t), \\
	A_0= -\sum_{k,l=1}^n a_{kl}(t,x)\partial_{x_k x_l}+\sum_{k=1}^n a_{k}(t,x)\partial_{x_k}+a_0,  
 \ A_i= \sum_{k=1}^n a_{k}^i(t,x)\partial_{x_k}+a_0^i, \ i=1,2,\ldots,r.
\end{multline*}
The equation   \eqref{e1} is furnished with 
the initial and boundary conditions 
\begin{equation}\label{e2}
	u|_{t=0}=u_0,\ \ Bu|_{S}=\frac{\partial u}{\partial \gamma}+\sigma u= g(t,x),\ \frac{\partial u}{\partial \gamma}= \sum_{i=1}^n \gamma_i(t,x) u_{x_i}.
\end{equation}
The additional conditions to determine the coefficients are as follows: 
\begin{equation}\label{e5}
	\int_{\Gamma}u(t,x)\varphi_j(x)=\psi_j(t),\    \  j=1,2,\ldots,s.
\end{equation}
The unknowns in the problem  \eqref{e1}-\eqref{e5} are a solution  $u$ and the functions 
$q_i(t)$ ($i=1,2,\ldots, s$).
The problems  \eqref{e1}-\eqref{e5} arise when describing
heat and mass transfer processes, diffusion and filtration processes, in ecology and many other
fields \cite{ali}, \cite{ozi}. The overdetermination conditions \eqref{e5} are not standard.   We can refer to the monograph
   \cite{pri} devoted to inverse parabolic problems and the monographs  \cite{bel,isa,kaba} containing basic statements of inverse problems including those for 
   parabolic equations.      In the monograph  \cite{bel} (see also the results in \cite{fro}), inverse problems of recovering  coefficients independent of some spatial variables are studied and 
         existence and uniqueness theorems are justified. In view of the method, the known  coefficients also are independent of some spatial variables. The overdetermination conditions are 
           values of a solution on some hyperplanes.  More complete results are exposed in  \cite{pya}-\cite{pya02}, where well-posedness of 
           inverse problems of recovering coefficients from the values of a solutions on some manifolds are studied. Inverse problems with pointwise or integral data 
           of the form  \eqref{e5} (in the latter case only the problem of identification of the heat transfer coefficient  $\sigma$ is considered in  \cite{kos}) are studied 
           in the articles by Prilepko A.I. and a series of interesting problems is described in \cite{pri}. Similar results but with other conditions on the data 
           and in other spaces are obtained in  \cite{pya1}-\cite{pya3}. Integral overdetermination conditions with  integrals  over the domain  $G$ are used in   \cite{py4}-\cite{py6}, where either the heat transfer coefficient or the flux on the boundary  are  to be determined.  These conditions is also often involved  in  \eqref{e5} \cite{kos},   \cite{dih,hao}, where some numerical methods for solving the problem are also exposed. 
In  \cite{lor} the problem of simultaneous recovering the heat transfer coefficient which depends on an integral of a solution with weight and lower-order coefficient  with the use of 
 the condition \eqref{e5} and the integral of a solution over the domain with weight are studied. 
Integral overdetermination conditions \eqref{e5} are often used as some approximations of pointwise conditions.  It is noted,  in particular, in  
  \cite{dih,hao}. Numerical solution of the problems  \eqref{e1}-\eqref{e5}  are also discussed in  \cite{sama}-\cite{hun1}). 
  The problem of recovering the heat transfer  coefficient with the overdetermination condition \eqref{e5} is considered in \cite{pya0}. Note that we do not know  articles devoted to the problem  \eqref{e1}-\eqref{e5} except for our
  article   \cite{pya00}. The main attention is paid to the case of recovering coefficients of higher order derivatives which are sought in the class of continuous functions. The existence and uniqueness theorems are obtained.
  But if we deal with real problems of heat and mass transfer then it is more interesting to consider the case of nonsmooth coefficients. 
  By the same  approach, we obtain here a similar result for the lower-order coefficients which belong to some Lebesgue space. Note that 
  the integral equations to which these two problems are reduced do not coincide. This

\section{Preliminaries}

Let $E$ be a Banach space. By  $L_p(G;E)$ ($G$
is a domain in  ${\mathbb R}^n$), we mean the space of measurable functions 
defined on  $G$ with values in 
$E$ and the finite norm  $\|\|u(x)\|_{E}\|_{L_p(G)}$ \cite{tri}.
We use  the Sobolev spaces   $W_p^s(G;E)$, $W_p^s(Q;E)$ 
 (see \cite{ama,ama1}) and the H\"{o}lder spaces $C^{\alpha,\beta}(\overline{Q})$, 
$C^{\alpha,\beta}(\overline{S})$ (see \cite{lad}). 
By a norm of a vector, we mean the sum of the norms of coordinates.  Given an interval $J=(0,T)$, put 
$W_p^{s,r}(Q)=W_p^{s}(J;L_p(G))\cap L_p(J;W_p^r(G)) $.
Respectively,  $W_p^{s,r}(S)=W_p^{s}(J;L_p(\Gamma))\cap
L_p(J;W_p^r(\Gamma))$.
 All function spaces and coefficients of the equation  \eqref{e1} are assumed to be real.
Next, we suppose that    $p>n+2$ and $\Gamma \in C^{3}$ (see the definition in \cite{lad}).
 Given  sets $X,Y$, the symbol $\rho(X,Y)$ stands for the distance between them.
Introduce the notations 
$Q^{\tau}=(0,\tau)\times G$,
$S^{\tau}=(0,\tau)\times \Gamma$. $G_\delta=\{x\in G:\rho(x,\Gamma)< \delta\}$, 
$Q_\delta=(0,T)\times G_\delta$, $Q_\delta^\tau=(0,\tau)\times G_\delta$.  
Construct a function  $\varphi(x)\in C^{\infty}(\overline{G})$ such that 
 $\varphi(x)=1$ in $G_{\delta/2}$ and $\varphi(x)=0$ in
$G\setminus  G_{3\delta/4}$. A small parameter  $\delta>0$ (it can be arbitrarily  small)  is fixed below. It characterizes a neighborhood of the boundary, 
where additional smoothness requirements are imposed on the data. 
Consider the  auxiliary problem 
	\begin{equation}\label{20}
		Mu =u_t +A_0u=f_0(t,x),\  \  A_0u= -\sum_{i,j=1}^n a_{ij}u_{x_i,x_j}+\sum_{i=1}^n a_iu_{x_i}+a_0u,
	\end{equation}
	\begin{equation}\label{18}
			\frac{\partial u}{\partial \gamma}+\sigma u\Bigl|_{S} =g, \ u|_{t=0} =u_0,
	\end{equation}
	We assume that 
\begin{equation}\label{6}
			u_0\in W_{p}^{2-2/p}(G),    B(0,x,D)u_0|_{\Gamma}=g(0,x),   g\in  W_{p}^{s_0,2s_0}(S), 	 f_0\in L_p(Q), 
			\end{equation}
where $s_0=1/2-1/2p$;
\begin{equation}\label{7}
	\varphi u_0(x) \in W_{p}^{3-2/p}(G),\ \varphi f_0\in L_p(0,T;W_p^1(G), \  g\in  W_{p}^{s_1,2s_1}(S),
\end{equation}
where $s_1=1-1/2p$;
 	\begin{equation}\label{eq33}
			a_{ij}\in 		C(\overline{Q}),\    \gamma_{i},\sigma \in W_p^{s_0,2s_0}(S),  a_{k}\in L_{p}(Q),     i,j=1,2,\ldots,n,\ k\leq n;
	\end{equation}
	\begin{equation}\label{351}
	a_{ij}\in L_{\infty}(0,T;W_{\infty}^{1}(G_{\delta})),\   a_{k}\in L_{\infty}(0,T;W_{p}^{1}(G_{\delta})),\  \gamma_{i},\sigma \in
		W_p^{s_1,2s_1}(S);
\end{equation}
where \ $i,j=1,\ldots,n$ \ k=0,1,\ldots,n,\  l=1,2,\ldots,r.
Next, we suppose  that 
\begin{equation}\label{e88}                                
|\sum_{i=1}^{n}\gamma_{i}\nu_{i}|\geq \varepsilon _{0}>0, \ \forall (t,x)\in S,
\end{equation}
where $\nu$ is the outward unit normal to  $\Gamma$ and $\varepsilon _{0}$ is a positive constant. 
The operator $A_0$ is elliptic, i.~e., there exists a constant $\delta_0>0$ such that 
\begin{equation}\label{e89}                                
\sum_{i,j=1}^n a_{ij}\xi_i\xi_j\geq \delta_0|\xi|^2\ \forall (t,x)\in Q, \ \forall \xi\in {\mathbb R}^n.
\end{equation}
The  next theorem is a consequence of Theorem  1 in  \cite{bar}.  
		
	\begin{theorem}\label{thr1} Let the conditions  {\rm \eqref{6}, \eqref{eq33}, \eqref{e88}, \eqref{e89}} hold.
		Then there exists a unique solution  $u\in W_p^{1,2}(Q)$ to the problem  {\rm \eqref{20}, \eqref{18}}. The following estimates is valid: 
		\begin{equation}\label{14}
			\|u\|_{W_{p}^{1,2}(Q)}\leq  c\bigl[\|u_0\|_{W_{p}^{2-2/p}(G)}+  \|f_0\|_{L_p(Q)}  +
\|g\|_{W_p^{s_0,2s_0}(S)}\bigr].
		\end{equation}
If additionally the conditions  \eqref{7},   \eqref{351} hold then a solution  is such that $\varphi u_t\in L_p(0,T;W_p^1(G))$, $\varphi u\in L_p(0,T;W_p^3(G))$. 
In the case of   $g=0$ and $u_0=0$ we have the estimates 
				\begin{equation}\label{15}
						\|u\|_{W_{p}^{1,2}(Q^\tau)}\leq
			c  \|f\|_{L_p(Q^\tau)},
		\end{equation}
\begin{multline}\label{151}
						\|u\|_{W_{p}^{1,2}(Q^\tau)} + \|\varphi u_t\|_{L_p(0,\tau;W_p^{1}(G))} + \|\varphi u\|_{L_p(0,\tau;W_p^{3}(G))}\leq \\
			c(\|f_0\|_{L_p(Q^\tau)}+\|\varphi f_0\|_{L_p(0,\tau;W_p^1(G))}),
		\end{multline}
		where the constant  $c$ is independent of  $f, \tau\in (0,T]$. 
  \end{theorem}

\section{Existence and uniqueness theorems}

Proceed with additional conditions on the data. We assume that  
\begin{equation}\label{8}
	\psi_{i}\in W_p^1(0,T),\  \psi_i(0)=\int_\Gamma  u_0\varphi_i\,d\Gamma, \ \varphi_i\in L_q(\Gamma), \  
    \frac{1}{q}+\frac{1}{p}=1,
 	      \end{equation}
 where $\ 	i=1,\ldots,s, $
 \begin{equation}\label{8a}
  \ f_m\in L_\infty(0,T;L_p(G))\cap L_\infty(0,T;W_p^1(G_\delta))\cap C([0,T];L_p(\Gamma)),\ 
\end{equation}
where $m=r+1,\ldots,s,$
\begin{equation}\label{8c} a_k^l \in L_{\infty}(0,T;W_{p}^{1}(G_{\delta}))\cap L_{\infty}(0,T;L_p(G))\cap  C([0,T];L_p(\Gamma)), 
\end{equation}
for all $ l=1,\ldots,r,\   k=0,1,\ldots,n. $
Assuming that the conditions of Theorem 1 and the above conditions are fulfilled, we can construct a solution to the problem \eqref{20}, \eqref{18}.  It possesses the properties  $\Phi\in W_p^{1,2}(Q)$, $\varphi \Phi_t\in L_p(0,T;W_p^1(G))$,
$\varphi \Phi\in L_p(0,T;W_p^3(G))$. 
Consider the matrix  $B_0$ of dimension  $s \times s$ with the rows 
\begin{multline*}
-<A_{1}(t,x,D) \Phi,\varphi_j>,\ldots,-<A_{r}(t,x,D)\Phi,\varphi_j>, \\
<f_{r+1}(t,x),\varphi_j>, \ldots, <f_{s}(t,x),\varphi_j>,\ j\leq s,
\end{multline*}
where $<u,v>=\int_{\Gamma}u(x)v(x)\,d\Gamma$. The embedding theorems imply that $\Phi\in C^{1-(n+2)/2p,2-(n+2)/p}(\overline{Q})$ (see Lemma 3.3, Ch.2 \cite{lad}).
In this case the entries of the matrix $B_0$ are continuous functions.  
We require that 
\begin{equation}\label{11}
	det\, B_0 \neq 0 \ \forall t\in [0,T].
\end{equation}
Since in Theorem \ref{thr2} below we prove existence of solutions locally in time, the condition \eqref{11} can be replaced with the condition
\begin{equation}\label{111}
	det\, B \neq 0,\  \forall t\in [0,T].
\end{equation}
where the matrix $B$ has the rows 
\begin{multline*}
-<A_{1}(t,x,D) u_0,\varphi_j>,\ldots,-<A_{r}(t,x,D)u_0,\varphi_j>, \\
<f_{r+1}(t,x),\varphi_j>, \ldots, <f_{s}(t,x),\varphi_j>,\ j\leq s.
\end{multline*}
This condition is easier to check in contrast to that in \eqref{11}. 
If we replace in Theorem \ref{thr2}  the condition \eqref{11} with the condition \eqref{111} then the claim of Theorem \ref{thr2} remains valid.

Make the change of variables $u=v+\Phi$. The problem \eqref{e1}-\eqref{e5}  is reduced to an equivalent problem
		\begin{equation}\label{28}
		v_t+A(\vec{q})v=\sum_{i=r+1}^s q_i f_i  -\sum_{i=1}^r q_i A_i \Phi=f^1,
	\end{equation}
where $A(\vec{q})v=A_0v+ A^1(\vec{q})v$,\ $A^1(\vec{q})=\sum_{i=1}^{r}q_iA_i$;
		\begin{equation}\label{28a}
		  v|_{t=0}=0, \  Bv|_S = 0;\  
	\end{equation}
\begin{equation}\label{32}
		<v,\varphi_j>=\psi_j(t)-<\Phi,\varphi_j>=\tilde{\psi}_j,\ \ i=1,\ldots,s.
	\end{equation}
					
\begin{theorem}\label{thr2} Let the conditions 
 {\rm \eqref{6}-\eqref{e89},   \eqref{8}-\eqref{11}} hold.  Then there exists a number  $\tau_0\in (0,T]$ such that 
 on the segment  $[0,\tau_0]$ there exists a unique solution   $(u,q_1,q_{2},\ldots,q_{s})$ to the problem 
		  {\rm \eqref{e1}-\eqref{e5}} such that 
		$u\in  W_{p}^{1,2}(Q^{\tau_0})$, $\varphi u\in L_{p}(0,\tau_0;W_{p}^{3}(G))$, $\varphi u_{t}\in
		L_{p}(0,\tau_0;W_{p}^{1}(G))$,   $q_j \in L_p(0,\tau_0)$, $ j=1,2,\ldots,s$.
\end{theorem}
	
	{\bf Proof.} 	We have reduced our problem  \eqref{e1}-\eqref{e5} to an equivalent simpler problem  \eqref{28}- \eqref{32}.
Denote by  $H_\tau$ the space of functions  $v$ satisfying the  conditions \eqref{28a} and such that 
$v\in W_{p}^{1,2}(Q^\tau),$ $\varphi v_t\in  L_p(0,\tau;W_p^{1}(G))$, $\varphi v\in L_p(0,\tau;W_p^{3}(G))$. Let 
$\|v\|_{H_\tau}=\|v\|_{W_{p}^{1,2}(Q^\tau)} + \|\varphi v_t\|_{L_p(0,T;W_p^{1}(G))} + \|\varphi v\|_{L_p(0,T;W_p^{3}(G))}$.
Define also the space    $W_\tau$  of functions  $f\in L_p(Q^\tau)$ such that  $\varphi f\in L_p(0,\tau;W_p^1(G))$.
Endow this space with the norm  $\|f\|_{W_\tau}=\|f\|_{L_p(Q^\tau)}+\|\varphi f\|_{L_p(0,\tau;W_p^1(G))}$. 
In view of Theorem  1, for every function  $f^1\in W_\tau$, there exists a unique solution 
$v=(\partial_t-A_0)^{-1}f ^1$ to the equation  $v_t-A_0v=f^1$ satisfying the conditions  \eqref{28a} and the estimate 
	\begin{equation}\label{332}
		\|v\|_{H_\tau}\leq 	c_1\|f^1\|_{W_\tau}, 
\end{equation}
where the constant  $c_1$ is independent of  $\tau$.
Reduce our problem to an integral equation. 
Multiply the equation \eqref{28} by $\varphi_j$ and integrate over $\Gamma$. Note that in the class of solutions $u$ described in Theorem 2  traces  of summands occurring in \eqref{28} on $\Gamma$ exist.
We obtain the equality 
		\begin{equation} \label{29}
		\tilde{\psi}_j'+ <A(\vec{q})v, \varphi_j>=  
		-\sum_{i=1}^{r}q_{i}<A_i\Phi,\varphi_j>
		+ \sum_{i=r+1}^{s}q_{i} <f_{i},\varphi_j>.
	\end{equation}
	The right-hand side of this equality coincides with  the $j$-th coordinate of the vector   $B_0(t) \vec{q}$.
	In this case the system  \eqref{29} can be written in the form 
	\begin{multline} \label{31}
		\vec{q}(t)=  B_0^{-1}H(\vec{q})(t)=R(\vec{q}), \   H(\vec{q})=(\tilde{\psi}_1'+<A(\vec{q}) v,\varphi_1>,
	\\	\tilde{\psi}_2' + <A(\vec{q}) v,\varphi_2>,\ldots, \tilde{\psi}_s'+<A(\vec{q}) v,\varphi_s>)^T,
	\end{multline}
 where $v$ is a solution to the problem \eqref{28}, \eqref{28a}. 
 Let  $\vec{q}=0$. In this case the right-hand side in \eqref{31} is written as follows:
 $$
 R(0)=B_0^{-1}\vec{\Psi},\  \vec{\Psi}=(\tilde{\psi}_1',\tilde{\psi}_2',\ldots, \tilde{\psi}_s')^T.
 $$
 Assign $R_0=2\|R(0)\|_{L_p(0,T)}$. Introduce that ball $B_{\tau}=\{\vec{q}\in L_p(0,\tau):\ \|\vec{q}\|_{L_p(0,\tau)}\leq R_0\}$.
 Next our aim is to prove that there exists a solution to the equation \eqref{31}. We employ the fixed point theorem.  
 To this end,
we first obtain estimates for solutions 	to the problem \eqref{28}, \eqref{28a} assuming that $\vec{q}\in B_\tau$. We have that
\begin{equation}\label{333}
v=-(\partial_t+A_0)^{-1} A^1(\vec{q})v+(\partial_t+A_0)^{-1}f^1.
	\end{equation}
From  \eqref{332}, it follows that 
\begin{equation}\label{15a}
		\|(\partial_t+A_0)^{-1} A^1(\vec{q})v\|_{H_\tau}\leq 	c_2\|A^1(\vec{q})v\|_{W_\tau}.
		\end{equation}
The following estimate is valid in accord with our condition on the coefficients:
\begin{multline} \label{335}
		\|A^1(\vec{q})v\|_{W_\tau} \leq 
		\|\sum_{j=1}^{r}q_{j} A_j(t,x,D)v\|_{W_\tau}		\\ \leq \|\vec{q}\|_{L_p(0,\tau)}\sum_{j=1}^{r}(\|A_j(t,x,D)v\|_{L_\infty(0,\tau;L_p(G))}+ 
 \|\varphi A_j(t,x,D)v\|_{L_\infty(0,\tau;W_p^1(G))}) \\ \leq  c_3\|\vec{q}\|_{L_p(0,\tau)}(\|v\|_{L_\infty(0,\tau;W_\infty^1(G))}+\|\varphi v\|_{L_\infty(0,\tau;W_\infty^2(G))}).  
			\end{multline}
The constant $c_3$ depends on the norm of coefficients in $L_\infty(0,\tau;L_p(G))\cap L_\infty(0,\tau;W_p^1(G_\delta))$.
Note that $W_p^{1,2}(Q^\tau)\subset C([0,\tau];W_p^{2-2/p}(G))$. Moreover, 
 if $\varphi v_t\in L_p(0,\tau;W_p^1(G))$, $\varphi v\in L_p(0,\tau;W_p^3(G))$ then $\varphi u\in C([0,\tau];W_p^{3-2/p}(G))$. These embeddings result from \cite[theorem III 4.10.2]{ama1}. 
  Even more,   in both cases the embedding constant is indepedent of $\tau$. 
Next, we have the estimate (see Theorem 4.6.1, 4.6.2 \cite{tri})
\begin{multline}\label{a1}
\|v\|_{L_\infty(0,\tau;W_\infty^1(G))}+ \|\varphi v\|_{L_\infty(0,\tau; W_\infty^{2}(G))}\leq c_4 
(\|v\|_{L_\infty(0,\tau;W_p^{1+s}(G))}+\\ \|\varphi v\|_{L_\infty(0,\tau; W_\infty^{2+s}(G))})\leq 
c_5(\|v\|_{L_\infty(0,\tau; W_p^{2-2/p}(G))}^\theta \|v\|_{L_\infty(0,\tau; L_p(G))}^{1-\theta} + \\ \|\varphi v\|_{L_\infty(0,\tau; W_p^{3-2/p}(G))}^{\theta_1} \|v\|_{L_\infty(0,\tau; L_p(G))}^{1-\theta_1}),
\end{multline}
where $n/p<s<1-2/p$, $1+s=\theta(2-2/p)$,  $2+s=\theta_1(3-2/p)$ and we have used the interpolation inequalities \cite{tri}.
The Newton-Leibnitz formula yields 
\begin{equation}\label{a2}
\|v\|_{L_\infty(0,\tau; L_p(G))}\leq \tau^{(p-1)/p}\|v_t\|_{L_p(Q^\tau)}.
\end{equation}
The  estimates \eqref{335}, \eqref{a1}, \eqref{a2} ensure the inequality  
\begin{equation} \label{a3}
  \|A^1(\vec{q})v\|_{W_\tau} \leq 
		 c_6\|\vec{q}\|_{L_p(0,\tau)}\|v\|_{H_\tau}\tau^{\beta},\  \beta=\min(\beta_0,\beta_1), 
  \end{equation}
  where   $\beta_0=(1-\theta)(p-1)/p$,   $\beta_1=(1-\theta_1)(p-1)/p$.
  Choose $\tau_0$ from the condition $c_6R_0\tau^{\beta}=1/2$.
  In this case the relation \eqref{333} yields   
  \begin{equation}\label{a4}
\|v\|_{H_\tau}\leq 2\|(\partial_t+A_0)^{-1}f^1\|_{H_\tau}.
	\end{equation}
In view of \eqref{332}, we can write out the estimate 
  \begin{equation}\label{a5}
\|v\|_{H_\tau}\leq c_7\|f^1\|_{W_\tau}\leq c_8\|\vec{q}\|_{L_p(0,\tau)},\ \tau\leq \tau_0,
	\end{equation}
where we can assume that the constant $c_8$ is independent of $\tau\leq \tau_0$. It depends on the norms 
of the data. 
  Next, we consider two vectors $\vec{q}_i=(q_1^i,q_2^i,\ldots,q_s^i)$, $i=1,2$, and construct the corresponding 
solutions $v_i$ to the problem     \eqref{28}, \eqref{28a}. 
  Subtracting two equations \eqref{28} we conclude that the difference $w=v_1-v_2$ is a solution to the problem 
  		\begin{multline}\label{a6}
		w_t+A((\vec{q}_1+\vec{q}_2)/2)w+ A^1(\vec{q}_1-\vec{q}_2)(v_1+v_2)/2=\\ \sum_{i=r+1}^s  (q_1^i-q_2^i) f_i  -\sum_{i=1}^r (q_1^i-q_2^i) A_i \Phi,
	\end{multline}
		\begin{equation}\label{a7}
		  w|_{t=0}=0, \  Bw|_S = 0.   
	\end{equation}
As before, in view of the inequality  \eqref{332}, we have the estimate
\begin{multline}\label{a8}
\|w\|_{H_\tau} \leq c_1\|A^1((\vec{q}_1+\vec{q}_2)/2)w\|_{W_\tau} + \\ c_1\|A^1(\vec{q}_1-\vec{q}_2)(v_1+v_2)/2\|_{W_\tau}+c_{10}\|\vec{q}_1-\vec{q}_2\|_{L_p(0,\tau)}. 
\end{multline}
Applying  the estimate \eqref{a3}, we infer 
\begin{multline}\label{a9}
\|w\|_{H_\tau} \leq c_6R_0\tau^{\beta}\|w\|_{H_\tau} + c_6\tau^{\beta} \|\vec{q}_1-\vec{q}_2\|_{L_p(0,\tau)}\|(v_1+v_2)/2\|_{H_\tau}+\\ c_{10}\|\vec{q}_1-\vec{q}_2\|_{L_p(0,\tau)}. 
\end{multline}
Since $\tau\leq \tau_0$, the last estimate implies that 
\begin{equation}\label{a9}
\|w\|_{H_\tau} \leq  2c_6\tau^{\beta} \|\vec{q}^1-\vec{q}_2\|_{L_p(0,\tau)}\|(v_1+v_2)/2\|_{H_\tau}+2c_{10}\|\vec{q}_1-\vec{q}_2\|_{L_p(0,\tau)}.
\end{equation}
Next, involving the estimate \eqref{a5} written for the functions $v_i$,
we conclude that 
      \begin{equation}\label{a10}
\|v_1-v_2\|_{H_\tau}\leq c_{11}\|\vec{q}_1-\vec{q}_2\|_{L_p(0,\tau)},
	\end{equation}
with a constant $c_{11}$ independent of $\tau\leq \tau_0$.
  Estimate the norm $\|R(\vec{q}_1)-R(\vec{q}_2)\|_{L_p(0,\tau)}$ with $\tau\leq \tau_0$. 
  Actually, the required estimate has been already established. 
  Consider the expressions $I_j=<(A(\vec{q}_1)v_1-A(\vec{q}_2))v_2,\varphi_j>$ occurring in the difference 
    $R(\vec{q}_1)-R(\vec{q}_2)$. 
    The H\"{o}lder inequality yields 
    \begin{multline*}
    |I_j|=|<A_0 w +A^1((\vec{q}_1+\vec{q}_2)/2)w+ A^1(\vec{q}_1-\vec{q}_2)(v_1+v_2)/2,\varphi_j>|\leq  \\
    c_{12}(\|A_0 w\|_{L_p(\Gamma)} +\|A^1((\vec{q}_1+\vec{q}_2)/2)w\|_{L_p(\Gamma)}+ \|A^1(\vec{q}_1-\vec{q}_2)(v_1+v_2)/2\|_{L_p(\Gamma)}.
    \end{multline*}
 The second summand on the right-hand side of $I_j$ is estimated as follows (see \eqref{a3}):
 \begin{multline}\label{a11}
 \|A^1((\vec{q}_1+\vec{q}_2)/2)w\|_{L_p(0,\tau;L_p(\Gamma))}\leq c_{12}\| A^1((\vec{q}_1+\vec{q}_2)/2)w\|_{W_\tau}\\ \leq 
 c_{13}\tau^\beta R_0\|w\|_{H_\tau}.
 \end{multline}
In view of the inequality \eqref{a5} written for the functions $v_i$, the third summand admits the estimate
 \begin{multline}\label{a12}
 \|A^1(\vec{q}_1-\vec{q}_2)\frac{(v_1+v_2)}{2}\|_{L_p(0,\tau;L_p(\Gamma))}\leq c_{14}\| A^1(\vec{q}_1-\vec{q}_2)\frac{(v_1+v_2)}{2}\|_{W_\tau}\leq \\
 c_{13}\tau^\beta C_1(R_0)\|\vec{q}_1-\vec{q}_2\|_{L_p(0,\tau)}.
 \end{multline}
 At last,  the first summand is estimated as follows: 
 \begin{multline}\label{a13}
\|A_0w\|_{L_p(0,\tau;L_p(\Gamma))}\leq c_{15}\|\varphi A_0w\|_{L_p(0,\tau;W_p^s(G))}\leq \\ c_{16}\|\varphi A_0w\|_{L_p(0,\tau; W_p^1(G))}^s \|\varphi A_0w\|_{L_p(Q^\tau)}^{1-s},
  \end{multline}
  where $ n/p<s<1.  $  
  Here we have used the interpolation inequality \cite{tri}. 
  Using the conditions on the coefficients, we obtain the estimate
  \begin{equation}\label{a14}
  \|\varphi A_0w\|_{L_p(G)}\leq c_{16}(\|\varphi w\|_{W_p^2(G)}+\| w\|_{L_\infty(0,\tau;W_\infty^1(G)}).
  \end{equation}
   Next, we repeat the arguments used in proof of the estimate \eqref{a3} (see the inequality \eqref{a1}).
Finally, we arrive at the estimate 
\begin{equation}\label{a15}
  \|\varphi A_0w\|_{L_p(G)}\leq c_{17}\tau^{\beta_2} \|\varphi w\|_{H_\tau}, \ \beta_2>0.
  \end{equation}
  The estimates \eqref{a11}, \eqref{a12}, \eqref{a15}, \eqref{a10} ensure the estimate 
  \begin{equation}\label{a16}
  \|R(\vec{q}_1)-R(\vec{q}_2)\|_{L_p(0,\tau)}\leq  c_{18}\tau^{\beta_3} \|\vec{q}_1-\vec{q}_2\|_{L_p(0,\tau)},
  \end{equation}
  where $\beta_3$ is a positive constant and $c_{18}$ is a constant depending on $R_0$ but independent of $\tau$. 
  Choose a quantity $\tau_1\leq \tau_0$ such that  $c_{18}\tau^{\beta_3}\leq 1/2$ for $\tau\leq \tau_1$. 
  In this case the operator $R(\vec{q})$ takes the ball $B_{\tau_1}$ into itself and is a contraction.
  The fixed point theorem implies that the equation \eqref{31} has a unique solution in the ball $B_{\tau_1}$. 
  Let   $v=v(\vec{q})$. Demonstrate that this functions satisfies 
\eqref{32}.
Integrating   \eqref{28} with the weight  $\varphi_j$ over $\Gamma$, we obtain the equalities 
\begin{multline} \label{42}
	<v_t,\varphi_j>+<A(\vec{q}) v,\varphi_j>=
	-\sum_{j=1}^{r}q_{j} <A_j\Phi,\varphi_j>  \\
	+ \sum_{j=r+1}^{s} <f_{j},\varphi_j>q_{j}(t).
\end{multline}
Subtracting them from  \eqref{29}, we conclude that  $<v_t,\varphi_j>-\tilde{\psi}_j'=0$ for all  $j$ and thereby the conditions  \eqref{32} hold.
Uniqueness of solutions follows from the estimates obtained in the above proof.

{\bf Remark. } Stability estimate for solutions also holds. 

In the linear case, i,e, the unknown functions occur only into the right-hand side the claim becomes	global in time.
In this case  $A=A_0$ and the matrix  $B_0$ has the rows  
\begin{equation}
<f_{1}(t,x),\varphi_j>, \ldots, <f_{s}(t,x),\varphi_j>,\ j\leq s,
\end{equation}
i.e., $r=0$.

More exactly, the following theorem is valid.

\begin{theorem} Let the conditions 
 {\rm \eqref{6}-\eqref{e89},   \eqref{8}, \eqref{8a}, \eqref{11}} hold.  Then there exists a unique 
  solution   $(u,q_1,q_{2},\ldots,q_{s})$ to the problem 
		  {\rm \eqref{e1}-\eqref{e5}} such that 
		$u\in  W_{p}^{1,2}(Q)$, $\varphi u\in L_{p}(0,T;W_{p}^{3}(G))$, $\varphi u_{t}\in
		L_{p}(0,T;W_{p}^{1}(G))$,   $q_j \in L_p(0,T)$, $ j=1,2,\ldots,s$.
\end{theorem}
	
The proof is  in  lines with that of the main result in \cite{pya02}.
The idea of the proof is as follows. First, we prove solvability of the problem   on some small segment  of time $[0,\tau_0]$. 
Next, repeating  the arguments we establish solvability on $[\tau_0,\tau_1]$, and so on.  Since the problem is linear, it is possible to prove that the length $\tau_i-\tau_{i-1}$
does not tends to zero which implies solvability on the whole segment $[0,T]$.

\section{Discussion}

Under consideration are mathematical models of heat and mass transfer.
We study  inverse problems of recovering lower-order  coefficients and the right-hand side in  a second order parabolic equation.
The coefficients  are representable in the form  of  a finite segments of series
  whose coefficients depending on time are to be determined. The linear case in which only the right-hand side is recovered  is also considered.
  The overdetermination conditions are the integrals over the boundary of a domain of a solution with weights. 
The main attention is paid to existence,  uniqueness, and stability estimates for
solutions to inverse problems of this type. The problem is reduced to an operator equation which
is studied with the use of the fixed point theorem and a priori estimates.   A solution has all generalized derivatives occurring into the equation summable to some power. 
 The method of the proof is constructive and it can be used for developing new numerical algorithms for solving the problem.

\section{Conclusions}

The existence and uniqueness theorems in inverse problems of recovering the lower-order coefficients in a parabolic equation  from the boundary integral measurements
are proven locally  in time. They  are sought in the form of a finite segment of the Fourier series with coefficients depending on time.  The proof relies on  a priori bounds and  the fixed point theorem.
The conditions on the data ensuring existence and uniqueness of solutions in Sobolev
classes  are  sharp. They are smoothness and consistency conditions on the
data and additional conditions on the kernels of the integral operators used in additional measurements.


\vspace{6pt}




Acknowledgement. {The research was carried out within the state assignment of Ministry of Science and Higher Education of the Russian Federation (theme No. FENG-2023-0004,
"Analytical and numerical study of inverse problems on recovering parameters of atmosphere or water pollution sources and (or) parameters of media").}

\vspace{5mm} \noindent {Sergey  Pyatkov, Doctor of Sci.,  Professor, Engineering School of  Digital Technologies, Yugra State University,  Khanty-Mansiysk, Russia, E-mail: pyatkovsg@gmail.com; \\
Oleg Soldatov, postraduate, Engineering School of  Digital Technologies, Yugra State University,  Khanty-Mansiysk, Russia, E-mail: oleg.soldatov.97@bk.ru } 


\begin{thebibliography}{999}


\bibitem {ali}
 Alifanov, O.M.; Artyukhin, E.A.; Nenarokomov, A.V.
{\it Inverse Problems in the Study of Complex Heat Transfer};
 Janus-K: Moscow, Russia, 2009.

\bibitem{ozi}  Ozisik,  M.N.;   Orlande, H.R.B.  {\it Inverse Heat Transfer};    Taylor \& Francis: New York, USA
        2000.

\bibitem{pri}
Prilepko, A.I.,  Orlovsky, D.G., Vasin, I.A. {\it Methods for Solving Inverse Problems in Mathematical Physics};   Marcel Dekker: New York,  1999.

\bibitem{bel}
Belov, Ya.Ya. {\it Inverse Problems for Parabolic Equations}; VSP:  Utrecht,
2002. 

\bibitem{isa}
Isakov, V. {\it  Inverse Problems for Partial Differential
Equations}; Springer: Berlin,  2006.

\bibitem{kaba}
Kabanikhin, S.I. {\it Inverse and Ill-Posed Problems}; Walter de Gruyter: Berlin, Boston,  2012. 

\bibitem{fro}   Frolenkov, I.V., Romanenko, G.V. 
    On Solving an Inverse Problem for a Multidimensional Parabolic Equation,
    {\it Siberian Journal of Industrial Mathematics},  {\bf 2012}, {\it 15}, no.\,2(50), 139--146. 

\bibitem{pya}  Pyatkov, S.G., Samkov, M.L.
    On Some Classes of Coefficient Inverse Problems for Parabolic 
    Systems of Equations,     {\it Sib. Adv. in Math.,} {\bf 2012}, {\it 22}(4), 287--302.

\bibitem{pya01}  Pyatkov, S.G., Tsybikov, B.N. 
    On Some Classes of Inverse Problems for Parabolic and Elliptic Equations
    {\it J. Evol. Equat.}, {\bf 2011}, {\it 11}(1), 155--186.

\bibitem{pya02}  Pyatkov, S.G. 
    On Some Classes of Inverse Problems for Parabolic Equations,
     {\it J. Inv.  Ill-Posed problems}, {\bf 2011}, {\it 18}(8), 917--934.

\bibitem{kos}
Kostin, A.B., Prilepko, A.I. On Some Problems of the Reconstruction of a
Boundary Condition for a Parabolic Equation, II,  {\it Differ. Equat.}, {\bf 1996}, {\it 32}, 1515--1525.


\bibitem{pya1}  Pyatkov, S.G. 
    On Some Classes of Inverse Problems with Overdetermination 
    Data on Spatial Manifolds,
    {\it Siberian Mathematical Journal,} {\bf 2016}, {\it 57}(5), 870--880.

\bibitem{pya2}  Pyatkov, S.G., Rotko, V.V. 
    On some parabolic inverse problems with the 
    pointwise overdetermination
    {\it Siber. Adv. in Math.}, {\bf 2020}, {\it 30}(2), 124--142.

\bibitem{pya4}   Pyatkov, S.G. 
    Identification of Thermophysical Parameters in Mathematical 
    Models of Heat and Mass Transfer,     {\it Journal of Computational and Engineering Mathematics,} {\bf 2022}, {\it 9}(2), 52--66.

\bibitem{pya3}  Pyatkov, S.G., Rotko, V.V. 
    Inverse Problems with Pointwise Overdetermination 
    for some Quasilinear Parabolic Systems,
    {\it AIP Conference Proceedings,} {\bf 2017}, {\it 1907}, 020008.

\bibitem{py4}
 Pyatkov, S.,  Soldatov, O.,  Fayazov, K.  Inverse problems of recovering the heat transfer coefficient with integral data, {\it  J. of Math. Sci.,} {\bf 2023}, {\it  274}(2), 255--268.
 
\bibitem{py5}
Kozhanov, A.I. Linear Inverse Problems for Some Classes of Nonlinear Nonstationary
Equations, {\it  Siberian Electronic Reports}, {\bf 2015}, {\it  12},  264--275. 

\bibitem{py6}
Verzhbitskiy,  M. A., Pyatkov, S. G. On Some Inverse Problems of Recovering Boundary Regimes, 
{\it Math. Notes of NEFU}, {\bf 2016}, {\it  23}(2),  3--18.

\bibitem{dih}
Dihn, N.; Hao, D.N.; Thanh, P.X.; Lesnik, D.
Determination of the Heat Transfer Coefficients in
Transient Heat Conduction. {\it  Inverse Problems,} {\bf 2013}, {\it 29},  095020.

\bibitem{hao}
Hao, D.N.;  Huong, B.V.; Thanh, P.X.; Lesnik, D.  Identification of Nonlinear
Heat Transfer Laws from Boundary Observations, {\it  Applicable Analysis,}  {\bf 2014}, {\it 94}(9), 1784--1799.

\bibitem{lor}
Lorenzi, A.; Messina, F.;  Identifying a Spherically Symmetric
Conductivity in a Nonlinear Parabolic Equation, {\it Applicable Analysis}, {\bf 85}(8), 867--889.

\bibitem{sama}  Samarskii, A.A.;  Vabishchevich P.N.
    {\it Numerical Methods for Solving Inverse Problems of Mathematical Physics.}
    Walter de Gruyter GmbH \& Co. KG: Berlin/Boston,   2007.

\bibitem{hun1}
 Huntul, M.J. Identification of the Timewise Thermal Conductivity in a 2D
Heat Equation from Local Heat Flux Conditions, {\it Inverse Problems in Science and Engineering}, 
{\bf 2021}, {it 29}(7), \,903--919. 

\bibitem{pya0}
	Pyatkov, S.G.; Soldatov, O.A.;  Identification of the Heat Transfer Coefficient from Boundary Integral Data, {\it Siberian Mathematical Journal,} 2024, {\bf  65}(4), 824--839.

\bibitem{pya00} Pyatkov, S.; Pronkina, T.  Coefficient Inverse Problems of Identification of Thermophysical Parameters from Boundary Integral Data, 
{\it Journal of Mathematical Sciences,} {\bf  282}(2),  241--252. 

\bibitem{tri}   Triebel H.
    {\it Interpolation Theory. Function Spaces. Differential Operators.}
    VEB Deutscher Verlag der Wissenschaften: Berlin, 1978.

\bibitem{ama}    Amann H. 
    Compact embeddings of vector-valued sobolev
	and besov spaces,     {\it Glasnik matematicki,} {\bf 2000}, {\it 35(55)}, 161--177.

\bibitem{ama1}    Amann H. 
    {\it Linear and Quasilinear Parabolic Problems},
    Birkhauser Verlag: Basel,   1995.

\bibitem{lad} Ladyzhenskaya, O.A.; Solonnikov, V.A.;   Uraltseva, N.N. {\it Linear and Quasilinear Equations of
    Parabolic Type}; American Mathematical Society: Providence, RI, USA,  1968.

\bibitem{bar}	Baranchuk, V.A.; Pyatkov, S.G. On some classes of inverse problems with pointwise ovedetermination for mathematical models of heat and mass transfer.
 Bulletin of the Yugra State University, {\bf 2020}, {\it 3}, 38--48.




\end{thebibliography}
\end{document}